\documentclass{amsart}

\newtheorem{theorem}{Theorem}[section]
\newtheorem{lemma}[theorem]{Lemma}
\newtheorem{corollary}[theorem]{Corollary}
\newtheorem{proposition}[theorem]{Proposition}

\theoremstyle{definition}
\newtheorem{definition}[theorem]{Definition}

\theoremstyle{remark}
\newtheorem{remark}[theorem]{Remark}

\newfont{\bbf}{msbm10 scaled\magstep1}

\def\Cc{\mbox{\bbf C}}

\def\ind{\mbox{\rm ind\,}}
\def\rank{\mbox{\rm rank\,}}

\def\deg{\mbox{\rm deg\,}}
\newcommand{\bc}[1]{\left(#1\right)}

\newcommand{\bs}[1]{\left\{ #1 \right\} }
 \newcommand{\bq}[1]{\left[#1\right]}
\newcommand{\m}[1]{\left|#1\right|}

\font\col=msbm10 scaled \magstep1

\newcommand{\C}{\mbox{\col C}}

\begin{document}

\title{On the Fredholm and Weyl spectrum of several commuting operators}
\author{R. Levy}
\address{Sofia University, Faculty of Mathematics and Informatics, Bd. J.Bourchier 5,
Sofia 1164, Bulgaria}
\email{LEVY@FMI.UNI-SOFIA.BG}


\thanks{The author was partially supported by the contract MM-1401/04 with the
Bulgarian Ministry of Science and Education.}


\date{}

\begin{abstract}
In the paper one considers the local structure of the Fredholm
joint spectrum of commuting $n$-tuples of operators. A connection
between the spatial characteristics of operators and the algebraic
invariant of the corresponding coherent sheaves is investigated. A
notion of Weyl joint spectrum of commuting $n$-tuple is
introduced.
\end{abstract}

\maketitle

\section{Introduction.}
Let $T=\left( T_1,\ldots,T_n \right)$, where $T_1,\ldots,T_n$ are
mutually commuting linear bounded operators acting in the Banach
space $X$. By $\sigma(T) \subset \Cc^n$ we will denote the joint
Taylor spectrum of $T$. Recall that $\sigma(T)$ consist of all
points $\lambda = \left( \lambda_1, \ldots ,\lambda_n \right)$ in
$\Cc^n$ such that the Koszul complex $K_*(T-\lambda,X)$ of the
operators $\left(T_1- \lambda_1, \ldots ,T_n-\lambda_n \right)$ is
not exact. Suppose that for given $\lambda \in \sigma(T)$ all the
homology spaces $H_i(T,\lambda):=H_i\bc{K_*(T-\lambda,X)}$ are
finite-dimensional: we will call such a point a {\it Fredholm}
point for $T$ and write $\lambda \in \sigma_F(T)$, the Fredholm
spectrum of $T$. The remaining part is called {\it essential
spectrum} and denoted by $\sigma_e(T):=\sigma(T) \backslash
\sigma_F(T)$.

It is a simple observation (see \cite{L})that any finite Fredholm
complex of Banach spaces with differentials, holomorphically
depending on the parameters, is holomorphically quasi-isomorphic
to a holomorphic complex of finite-dimensional spaces. Therefore,
the sheaves $\mathcal H_i(T)$ of homologies of the complex of
germs of holomorphic functions with values in $K_*(T-\lambda,X)$
are coherent on $\Cc^n \backslash \sigma_e(T)$. Fixing a point
$\lambda^0 \in \sigma_F(T)$, one can consider the stalk of the
homology sheaf $\mathcal H_i(T)_{\lambda^0}$ as a module over the
Noetherian local ring $\mathcal O_{\lambda^0}$ of germs of
holomorphic functions at this point. The Fredholm spectrum
$\sigma_F(T)$ is a complex-analytic subspace of $\Cc^n\backslash
\sigma_e(T)$, its dimension near $\lambda^0$ is well defined, and
is an integer not exceeding the dimension of the ambient space
$\Cc^n$; in the present paper we are interested mostly in the case
when this dimension is strictly less than $n$.

To any coherent sheaf, or finitely generated $\mathcal
O_{\lambda^0}$-module, one can attach an element of the cycle
group of $\mathcal O_{\lambda^0}$, i.e. a formal sum of prime
ideals of $\mathcal O_{\lambda^0}$.Taking the alternated sum of
cycles of the modules $\mathcal H_i(T)_{\lambda^0}$ for
$i=0,\ldots,n$, one obtains the cycle of the Koszul complex of the
$n$-tuple. Roughly speaking, any irreducible component of
$\sigma_F(T)$ through $\lambda^0$ "participates" in the homology
sheaves multiplied by some integer coefficient, and we obtain a
set of integers characterizing the homology sheaves of the Koszul
complex of $T$. (In the case of a single operator $T$ the
corresponding invariant is the index of $T-\lambda^0$.) One of the
purposes of the paper is to establish some connections between the
these algebraic characteristics of the homology sheaves, and the
action of operators of $T$ in $X$.

In the first section of the paper we recall some useful facts from
the commutative algebra and the theory of analytic local rings.
Especially, we recall the notions of the isolated prime ideals in
the support of a module, the cycle of a module, and the
functoriality of the cycle map under finite morphisms. We
introduce also the notions of the Hilbert-Samuel polynomial and
the multiplicity of a module with respect to a prime ideal.

The second section contains the main results of the paper. We give
a definition of the cycle $z_{\lambda^0}\bc{T}$ in a point
$\lambda^0$ of the Fredholm spectrum of $T$; this is a formal
linear combination with integer coefficients of the irreducible
components of $\sigma_F(T)$ passing through $\lambda^0$. In other
words, to any irreducible component $\mathcal P$ of $\sigma_F(T)$
containing $\lambda^0$ with corresponding prime ideal $p_{\mathcal
P}$ one attaches an integer $l_{\mathcal P}(T)$ with
$z_{\lambda^0}\bc{T}= \sum l_{\mathcal P}(T) \ p_{\mathcal P}$.
Theorem 2.3 proves the functoriality of the maximal component of
the cycle under the holomorphic functional calculus. Next,
theorems 2.5 and 2.8 establish some connections between the
coefficients of the cycle $z_{\lambda^0}\bc{T}$ and the spacial
characteristics of the $n$-tuple $T$; in some particular cases
this can be used in the calculations of these coefficients. In
particular, theorem 2.8, considering the structure of the last
homology sheaf of the Koszul complex of $T$, can be applied to the
theory of multidimensional contractions, developed by Arveson
(\cite{A1}, \cite{A2}); in more details the connection between the
characteristics of a multidimensional contraction and the geometry
of its Fredholm spectrum will be considered elsewhere.

The last part of the second section is devoted to the Weyl
spectrum of the commuting $n$-tuple $T$ of operators. In the paper
\cite{P}, M. Putinar defines the Weyl spectrum $\omega(T)$ of $T$
as the complement in $\sigma(T)$ of the set of points $\lambda \in
\sigma_F(T)$ such that $\ind (T-\lambda) = 0$. However, in the
case when $\dim \sigma_F(T)<n$, it is easy to see that $\ind
(T-\lambda) = 0$ for any $\lambda \notin \sigma_e(T)$, and
therefore the Weyl spectrum will coincide with the essential
spectrum  of the $n$-tuple. We propose here an alternative
definition: a Fredholm point $\lambda$ is not in the Weyl spectrum
iff the cycle $z_{\lambda}\bc{T}$ is zero, i.e. all $l_{\mathcal
P}(T)=0$. In the case when all the components of $\sigma_F(T)$ are
of (maximal) dimension $n$, this definition coincides with the
definition adopted in \cite{P}. However, even in the case of a
single operator there is a certain discrepancy between our
approach and the standard one: in our definition, all the isolated
point of the Fredholm spectrum of $T$ belong to its Weyl spectrum.
Further we give some basic properties of the so-defined Weyl
spectrum: any $n$-tuple with SVEP has the "Weyl property". We
prove also the property of spectral inclusion under holomorphic
functional calculus for this spectrum. In some particular cases
one has a stronger result - the spectral mapping theorem (see
propositions 2.18 and 2.19).

The present paper continues the investigations on the structure of
the Fredholm spectrum of a commuting $n$-tuple of operators
started in \cite{L}, and systematically uses some basic result of
that paper.

\section{Some necessary facts from the commutative algebra.}

Consider a coherent analytic sheaf $\mathcal L$, defined on open
subset of $\Cc^n$. Recall that its (geometric) support $\mbox{\rm
supp} \, {\mathcal L}$ consists of all point where $\mathcal L$ is
non-zero; $\mbox{\rm supp} \, {\mathcal L}$ is complex-analytic
subset of $\Cc^n$. If $\lambda^0 \in \mbox{\rm supp} \, {\mathcal
L}$, then the stalk of $\mathcal L$ at $\lambda^0$ is a finitely
generated module over the local Noetherian ring ${\mathcal
O}_{\lambda^0}$. In this section we recall for a latter use some
basic notions from the commutative local algebra, concerning the
Noeterian rings and modules.

Let $\mathcal A$ be a commutative local Noetherian ring, and
$\mathcal M$ be finitely generated $\mathcal A$-module with
annihilator $Ann({\mathcal M})$. Let $Prime\bc{\mathcal A}$ be the
set of prime ideals of $\mathcal A$. One denotes by $Supp
\bc{{\mathcal M}}$ the (algebraic) support of $\mathcal M$ - the
set of all prime ideals in $\mathcal A$, containing $Ann({\mathcal
M})$. The set $Supp({\mathcal M})$ has an natural ordering by
inclusion and the minimal elements with respect of this ordering
(called isolated associated primes) play a special role and
correspond to irreducible components of the geometric support of
${\mathcal M}$. The set $Iso({\mathcal M})$ of all such prime
ideals is finite (see \cite{E}, chap. 3).

For $p \in Prim\bc{\mathcal A}$, denote by ${\mathcal M}_p$ the
localization of the $\mathcal A$-module $\mathcal M$ with respect
to $p$. Then $p \in Supp \bc{{\mathcal M}}$ if and only if
${\mathcal M}_p \neq 0$. Moreover, $p \in Iso \bc{{\mathcal M}}$
iff $l_p\bc{\mathcal M} = \dim_{{\mathcal A}_p}\ {\mathcal M}_p$
is finite and nonzero. Suppose that $0={\mathcal M}_0 \subset
\ldots \subset {\mathcal M}_k = {\mathcal M}$ is a composition
series for ${\mathcal M}$, i.e. for all $i$ one has ${\mathcal
M}_i/{\mathcal M}_{i-1}= {\mathcal A}/p_i$ with $p_i \in
Prim\bc{{\mathcal A}}$. If $p \in Iso\bc{{\mathcal M}}$, then $p$
appears exactly $l_p\bc{{\mathcal M}}$ times in the sequence
$\bs{p_1,\ldots,p_k}$. By $Iso_{max}\bc{{\mathcal M}}$ we will
denote the set of elements of $Iso\bc{{\mathcal M}}$ of maximal
dimension (equal to the dimension of $\mathcal{M}$).

Let $\mathcal{M}$ and $\mathcal{A}$ be as above, and let $q$ be an
ideal in $\mathcal{A}$ such that $\mathcal{M}/q \mathcal{M}$ is of
finite dimension. Then there exists a polynomial
$P_{q,\mathcal{M}}$ (called the {\em Hilbert - Samuel polynomial}
of the pair $q,\mathcal{M}$) of degree $r=\dim_{\mathcal{A}}
\mathcal{M}$ such that for $n$ sufficiently big one has $\dim
\bc{\mathcal{M}/q^n \mathcal{M}}=P_{q,\mathcal{M}}(n)$. The
leading term of $P_{q,\mathcal{M}}$ has the form
$\frac{e_q\bc{\mathcal{M}}}{r!}n^r$, where $e_q\bc{\mathcal{M}}$
is a positive integer, called a {\em multiplicity} of $q$ at
$\mathcal{M}$. If $q=\mathfrak{m}$, where $\mathfrak{m}$ is the
maximal ideal of $\mathcal{A}$, then the corresponding integer is
called {\em multiplicity} of $\mathcal{M}$ and is denoted by
$e\bc{\mathcal{M}}$. The formula $\bc{*}$ in the page 1 of the
book \cite{S} shows how $e_q\bc{\mathcal{M}}$ can be calculated in
homological terms.

Suppose now that
$\overline{\mathcal{M}}=\bs{\mathcal{M}_0,\ldots,\mathcal{M}_n}$
is an ordered finite set of finitely generated modules over
$\mathcal{A}$. Denote $Ann\bc{\overline{\mathcal{M}}} =
\bigcap_{i=0}^n Ann\bc{\mathcal{M}}$. Let
$Supp\bc{\overline{\mathcal{M}}}$ be the set of all prime ideals
of $\mathcal{A}$, containing $Ann\bc{\overline{\mathcal{M}}}$,
$Iso\bc{\overline{\mathcal{M}}}$ -- the set of minimal elements of
$Supp\bc{\overline{\mathcal{M}}}$, and
$Iso_{max}\bc{\overline{\mathcal{M}}}$ -- the set of elements of
$Iso\bc{\overline{\mathcal{M}}}$ of maximal dimension. Take $i \in
\bs{0,\ldots,n}$ and a prime ideal $p \in Supp\bc{\mathcal{M}_i}$.
Since $Supp\bc{\mathcal{M}_i}\subset
Supp\bc{\overline{\mathcal{M}}}$, then if $p$ is not minimal in
$Supp\bc{\mathcal{M}_i}$, it will be not minimal in
$Supp\bc{\overline{\mathcal{M}}}$ as well. So, one obtains

\begin{proposition}
For any $p \in Iso\bc{\overline{\mathcal{M}}}$ and $i \in
\bs{0,\ldots,n}$ one has one of the following two possibilities:
1/ $p \in Iso\bc{\mathcal{M}_i}$, or 2/ $p \notin
Supp\bc{\mathcal{M}_i}$.
\end{proposition}

Therefore for any prime $p \in Iso\bc{\overline{\mathcal{M}}}$ and
$i$ the integer $l_p\bc{{\mathcal M}_i} = \dim_{{\mathcal A}_p}\
\bc{{\mathcal M}_i}_p$ is well-defined (it is positive in the case
1/ and zero in the case 2/). So for such a $p$ one can define:

$$l_p\bc{\overline{{\mathcal M}}}\ := \ \sum_{i=0}^n \ (-1)^i \ l_p\bc{{\mathcal
M}_i}.$$

Let $q$ be an ideal in $\mathcal{A}$ such that $\mathcal{M}_i/q
\mathcal{M}_i$ is finite-dimensional for all $i$. One can
introduce the Hilbert-Samuel polynomial of $\overline{M}$:

$$ P_{q,\overline{\mathcal{M}}}(n)\ := \ \sum_{i=0}^n \ (-1)^i \
P_{q,\mathcal{M}}(n).$$

If one denote by $r=\dim_{\mathcal{A}}\overline{\mathcal{M}}$ the
maximal dimension of the $\mathcal{A}$-modules $\mathcal{M}_i$,
$i=0,\ldots,n$, then $ P_{q,\overline{\mathcal{M}}}(n)$ is a
polynomial of dimension $r$ and with leading term
$\frac{e_q\bc{\overline{\mathcal{M}}}}{r!}n^r$, where

$$e_q\bc{\overline{\mathcal{M}}} \ := \ \sum_{i=0}^n \ (-1)^i \
e_q\bc{\mathcal{M}_i}.$$

Denote by $Z\bc{\mathcal A}$ the group of the cycles of $\mathcal
A$, i.e. of all formal linear combinations of elements of
$Prime\bc{\mathcal A}$ with integer coefficients. Then to any
tuple of $\mathcal{A}$-modules $\overline{\mathcal{M}}$ one can
attach the elements $z\bc{\overline{\mathcal{M}}}$ and
$z^{max}\bc{\overline{\mathcal{M}}}$ of $Z\bc{\mathcal A}$,
defined by the formulae

$$z\bc{\overline{\mathcal{M}}} \ = \ \sum_{p \in Iso
\bc{\overline{\mathcal{M}}}} l_p\bc{\overline{\mathcal {M}}}\
p\quad,\quad z^{max}\bc{\overline{\mathcal{M}}} \ = \ \sum_{p \in
Iso_{max} \bc{\overline{\mathcal{M}}}} l_p\bc{\overline{\mathcal
{M}}}\ p$$ and one has $$e_q\bc{\overline{\mathcal{M}}}\ = \
\sum_{p \in Iso_{max} \bc{\overline{\mathcal{M}}}} \
l_p\bc{\overline{\mathcal {M}}} e_q\bc{\mathcal A/p}$$ (see
\cite{S}, p. 125-126).

Finally, we will need a functoriality result. Let ${\mathcal B}$
be a subring of the analytic local ring ${\mathcal A}$ such that
the monomorphism $\varphi: {\mathcal B}\to {\mathcal A}$ is finite
in the sense of \cite{G-R}, II.2.2, i.e. ${\mathcal A}$ is
finitely generated as a ${\mathcal B}$-module. For any finitely
generated ${\mathcal A}$-module ${\mathcal M}$, denote by
${\mathcal M}_{\mathcal B}$ the underlying ${\mathcal B}$-module.
Then ${\mathcal M}_{\mathcal B}$ is a finitely generated
${\mathcal B}$-module again. For any $p \in Prime\bc{{\mathcal
A}}$, denote by $\varphi_*p$ the element $z^{max}\bc{p_{\mathcal
B}}$ of $Z\bc{\mathcal B}$. By linearity one can extend this
mapping up to a morphism $\varphi_* :Z\bc{\mathcal A} \to
Z\bc{\mathcal B}$.

\begin{proposition}

For any finitely generated ${\mathcal A}$-module ${\mathcal M}$
one has $z^{max}\bc{{\mathcal M}_{\mathcal
B}}=\varphi_*z^{max}\bc{{\mathcal M}}$.

\end{proposition}

\begin{proof}
Indeed, theorem 2 of \cite{G-R}, II.5.1 asserts that
$\dim_\mathcal A \mathcal L = \dim_\mathcal B \mathcal L_\mathcal
B$ for any finitely generated ${\mathcal A}$-module ${\mathcal
L}$. In particular, this is true for the modules of the type
${\mathcal A}/p$, $p \in Prime\bc{\mathcal A}$. Taking the
composition series for ${\mathcal M}$ and using the additivity of
the mapping $z^{max}$ (see \cite{S}, p. 125), we obtain
$z^{max}\bc{\mathcal{M}_\mathcal B}  =  \sum_{p \in Iso_{max}
\bc{\mathcal{M}}} l_p\bc{\mathcal {M}}\ p_\mathcal B =
\varphi_*z^{max}\bc{{\mathcal M}}$.
\end{proof}

\section{Local structure of the Fredholm and Weyl spectra of $n$-tuple}

\subsection{Main definitions}Consider a commuting $n$-tuple $T = \left(T_1,\ldots,T_n
\right)$ of operators acting in the Banach space $X$, and denote
by $K_*(T-\lambda,X)=
\left\{X_k,d_k(\lambda)\right\}_{k=0,1,\ldots,n} $ the Koszul
complex of $T$ in the point $\lambda \in \C^n$. Recall that $X_k$
is a direct sum of $n \choose k$ copies of the space $X$, and
$d_k(\lambda)$ depend linearly on $\lambda$. Denote by $\mathcal O
X$ the sheaf of germs of $X$-valued holomorphic functions on
$\C^n$; then one can consider on $ \C^n$ the complex $\mathcal O
K_*(T-\lambda,X)$ of sheaves of holomorphic sections of the
complex $K_*(T-\lambda,X)$. Denote by $H_i(T,\lambda)$ the $i$-th
homology space of the complex $K_*(T-\lambda,X)$, and by $\mathcal
H_i(T)$ the $i$-th sheaf of homologies of the complex of sheaves
$\mathcal O K_*(T-\lambda,X)$. Let $\mathcal H_i(T)_{\lambda^0}$
be the stalk of $\mathcal H_i(T)$ at ${\lambda^0}$, considered as
a local ${\mathcal O}_{\lambda^0}$-module.

Since the operators $T_i$ commute with the differentials
$d_k(\lambda)$ of the parameterized Koszul complex, one can define
the action of operators $T_k, k=1, \ldots , n $ on the spaces
$H_i(T,\lambda^0)$ and $\mathcal H_i(T)$. It is easy to see that
the operator $T_k$ acts on the space $H_i(T,\lambda^0)$ as a
multiplication by $\lambda^0_k$, and on the sheaf $\mathcal
H_i(T)$ - as a multiplication by the variable $\lambda_k$.

Fix a point $\lambda^0 \in \sigma_F(T)$; this means that all the
spaces $H_i(T,\lambda^0)$ are finite-dimensional. As it was noted
in \cite{L}, there exists a holomorphic complex of
finite-dimensional spaces, defined near $\lambda^0$ and
holomorphically quasi-isomorphic to $K_*(T-\lambda,X)$. Therefore,
all the sheaves $\mathcal H_i(T)$ are coherent in a neighborhood
of $\lambda^0$, and the stalks $\mathcal H_i(T)_{\lambda^0}$ at
$\lambda^0$ are finitely generated modules over the local
Noetherian ring $\mathcal O_{\lambda^0}$. Denote by
$\overline{\mathcal{H}}(T)_{\lambda^0}$ the $n+1$-tuple
$\bs{\mathcal H_i(T)_{\lambda^0}}_{i=0,\ldots,n}$. Then the prime
ideals in $Iso\bc{\overline{\mathcal{H}}(T)_{\lambda^0}}$
correspond to the irreducible components of the complex set
$\sigma_F(T)$, containing the point $\lambda^0$. Denote by $r=r
\bc{\lambda^0}$ the maximal dimension of these components.

\begin{definition}
Define
$$z_{\lambda^0}\bc{T}=z\bc{\overline{\mathcal{H}}(T)_{\lambda^0}}\,
, \,z_{\lambda^0}^{max}\bc{T}=
z^{max}\bc{\overline{\mathcal{H}}(T)_{\lambda^0}}.$$

\end{definition}

The element $z_{\lambda^0}\bc{T}$ has the form
$$z_{\lambda^0}\bc{T}\ = \ \sum_{p \in Iso
\bc{\overline{{\mathcal{H}}}(T)_{\lambda^0}}} l_p\bc{T}\  p$$
where $l_p\bc{T}=l_p\bc{\overline{{\mathcal{H}}}(T)_{\lambda^0}}$
are integers. We will call the integer $l_p\bc{T}$ a {\em local
index} of the $n$-tuple $T$ at the point $\lambda^0$ and the prime
ideal $p$.

This definition can be considered from the geometric point of
view. Let $\mathcal P$ be an irreducible component of the complex
set $\sigma_F(T)$; then in any point $\lambda \in \mathcal P$ it
determines a prime ideal $p_{\lambda}$ in the ring $\mathcal
O_\lambda$.

\begin{proposition}
The integer $l_{p_{\lambda}}(T)$ does not depend on the choice of
the point $\lambda \in \mathcal P$.
\end{proposition}

\begin{proof}
It is sufficient to proof that the integers $l_{p_{\lambda}}(T)$
are locally constant. Choose $\lambda^0 \in \mathcal P$ and a
sufficiently small ball $U$, centered at $\lambda^0$, such that
$U$ does not intersect any irreducible component of $\sigma_F(T)$,
not containing $\lambda^0$. Then for any $i \in \bs{0,\ldots,n}$
the $\mathcal O_U$-module $\mathcal H_i(T)_U$ has the same set of
isolated primes, and the same composition series, as $\mathcal
H_i(T)_{\lambda^0}$. The irreducible complex set $\mathcal P$
determines a prime ideal $p_U$ in $\mathcal O_U$. For any $\lambda
\in U$ the $\mathcal O_\lambda$-module $\mathcal
H_i(T)_{\lambda}$. If $\lambda \in \mathcal P$, the localization
of $p_U$ at $\lambda$ is a nontrivial prime ideal, and the
invariance of the composition series under localization shows that
$l_{p_{\lambda^0}}(T) = l_{p_{U}}(T)= l_{p_{\lambda}}(T)$.
\end{proof}

Now, to any irreducible component one can attach its local index
$l_\mathcal P(T)$. The set of all irreducible components and the
corresponding local indexes contains certain information about the
homology sheaves of the Koszul complex of $T$ and will be called a
{\em spectral picture} of the commuting $n$-tuple $T$. In the rest
of paper we will turn back to the algebraic point of view on the
local indexes and will give some facts allowing to compute it in
some cases.

\subsection{Particular cases: dimensions $1$ and $n$}
Suppose that the dimension of $\sigma_F(T)$ at $\lambda^0$ is the
maximal one, i.e. equal to $n$. This means that the Fredholm
spectrum $\sigma_F(T)$ contains a neighborhood $U$ of $\lambda^0$.
It is easy to see (or to derive from theorem 2.5 below) that in
this case
$$z\bc{T}_{\lambda^0}=z^{max}\bc{T}_{\lambda^0}=\bc{-1}^n\ind
\bc{T-\lambda^0} \ . \ \bq{0},$$ where $\ind \bc{T-\lambda^0}$ is
the Euler characteristic of the complex $K_*\bc{T-\lambda^0,X}$
and $\bq{0}$ is the zero ideal in ${\mathcal O}_{\lambda^0}$.

Suppose in addition that the $n$-tuple $T$ has SVEP. Then from the
theory of coherent sheaves it follows that the sheaf $\mathcal
H_n(T)$ is a free $\mathcal O$-module on the open dense subset
$\widetilde{U}$ of $U$. (In fact, $\widetilde{U}$ is a complement
of a complex subset of $U$ of dimension $<n$.) Then for any
$\lambda = \bc{\lambda_1, \ldots , \lambda_n} \in \widetilde{U}$
one has
$$\ind \bc{T-\lambda^0} = \ind \bc{T-\lambda} = (-1)^n \dim
\bc{X/\sum_{i=1}^n\bc{T_i-\lambda_i}X}.$$

Now consider the case when $\sigma_F(T)$ is one-dimensional near
$\lambda^0$. Then, as proven in \cite{L}, the modules $\mathcal
H_i(T)_{\lambda^0}$ are nonzero only for $i=n,n-1$. Moreover,
prop. 3.6 of \cite{L} shows that $$\mathcal H_{n-1}(T)\ \sim \
\mathcal H_n\bc{T^*},$$ where $T^*$ is the $n$-tuple
$\bc{T_1^*,\ldots,T_n^*}$ acting in the dual space $X^*$, and
$\sim$ denotes the equivalence modulo sheaves with
zero-dimensional support, i.e. concentrated in the point
$\lambda^0$. So one obtains the equality
$$z^{max}\bc{T}_{\lambda^0}\ =\ (-1)^n\bc{z^{max}\bc{\mathcal H_n(T)_{\lambda^0}}\ -
\ z^{max}\bc{\mathcal H_n\bc{T^*}_{\lambda^0}}}.$$

Finally, note that if $\sigma_F(T)$ is of dimension zero at
$\lambda^0$ (i.e. $\lambda^0$ is an isolated point of the Fredholm
spectrum), then $l_{\bs{\lambda^0}}(T)$ coincides with the
dimension of the spectral subspace of $T$ corresponding to the
point $\lambda^0$.

\subsection{Functoriality} Let
$f\bc{\lambda}=\bc{f_1\bc{\lambda},\ldots,f_k\bc{\lambda}}$ be a
$k$-tuple of holomorphic functions, defined in a neighborhood of
the spectrum $\sigma(T) \subset \Cc^n$ of the commuting $n$-tuple
$T$ of operators. It is well-known that
$\sigma(f(T))=f\bc{\sigma(T)}$ and
$\sigma_e(f(T))=f\bc{\sigma_e(T)}$. Take the point $\mu^0 \in
\Cc^k$ belonging to the Fredholm spectrum $\sigma_F(f(T))$ of the
operator $k$-tuple $f(T)$. Then its preimage $f^{-1}\bc{\mu^0}\cap
\sigma(T)$ is a finite subset $\bs{\lambda_1,\ldots,\lambda_p}$ of
the Fredholm spectrum $\sigma_F(T)$ of $T$. For any $\lambda \in
f^{-1}\bc{\mu^0}\cap \sigma(T)$ the inverse image $f^*$ by $f$
provides an embedding ${\mathcal O}^k_{\mu^0}\to {\mathcal
O}^n_{\lambda}$, and the induced monomorphism ${\mathcal
O}^k_{\mu^0}\to {\mathcal
O}^n_{\lambda}/Ann\bc{\overline{\mathcal{H}}(T)_{\lambda}}$ is
finite. Denote by $f_*: Z\bc{{\mathcal
O}^n_{\lambda}/Ann\bc{\overline{\mathcal{H}}(T)_{\lambda}}} \to
Z\bc{{\mathcal O}^k_{\mu^0}}$ the morphism used in prop. 1.2.

\begin{theorem}
One has $$ z^{max}\bc{\overline{\mathcal{H}}(f(T))_{\mu^0}}\ = \
\sum_{\lambda \in f^{-1}\bc{\mu^0}\cap \sigma(T)} f_*
z^{max}\bc{\overline{\mathcal{H}}(T)_{\lambda}}$$
\end{theorem}

\begin{proof}
In proposition 3.2 of \cite{L}, one states the equality
$${\mathcal{H}_i}(f(T))_{\mu^0}= \bigoplus_{\lambda \in f^{-1}\bc{\mu^0}\cap
\sigma(T)}{\mathcal{H}_{i+n-k}}(T)_\lambda \ ,$$ both sides
considered as ${\mathcal O}^k_{\mu^0}$-modules. Applying the
functor $z^{max}$ and taking into account the equality
$z^{max}\bc{{\mathcal M}_{\mathcal
B}}=\varphi_*z^{max}\bc{{\mathcal M}}$ proved in the previous
section, one obtains the theorem.
\end{proof}

\subsection{Computations for $z^{max}_{\lambda^0}\bc{T}$}
We will try to extract some information on the cycle
$z^{max}\bc{T}_{\lambda^0}$, i.e. on the integers
$l_p(T):=l_p\bc{\overline{\mathcal{H}}(T)_{\lambda^0}}$ for $p \in
Iso_{max}\bc{\overline{\mathcal{H}}(T)_{\lambda^0}}$, from the
properties of operators of $T$. Suppose that the coordinates
$\lambda' = \bc{\lambda_1,\ldots,\lambda_r}$ form a coordinate
system for $\sigma_F(T)$ at $\lambda^0$. This means that
$\lambda^0$ is an isolated point of $\Pi^{-1}\lambda'^0\cap
\sigma_F(T)$, where $\Pi$ is the coordinate projection onto first
$r$ coordinates, and
$\lambda'^0=\bc{\lambda^0_1,\ldots,\lambda^0_r}$,
$\lambda''^0=\bc{\lambda^0_{r+1},\ldots,\lambda^0_n}$. This always
can be achieved by a small perturbation of the coordinate system.
We will denote $T=\bc{T',T''}$, where $T'=\bc{T_1,\ldots,T_r}$,
$T''=\bc{T_{r+1},\ldots,T_n}$.

Consider the Koszul complex $K_*(T'-\lambda'^0,X)$ of the
operators $T_1-\lambda^0_1,\ldots,T_r-\lambda^0_r$ in $X$. Its
homology spaces $H_i(T',\lambda'^0)$, $i=0,\ldots,r$ are in
general not separated, but because of the commutativity of $T$ the
action of the operators from $T''$ on it is correctly defined.
Moreover, one can form the parameterized Koszul complex of
operators of $T''-\lambda''$ in $H_i(T',\lambda'^0)$; this complex
is exact for $\lambda'' \notin \Pi^{-1}\lambda'^0\bigcap
\sigma(T)$. The following proposition is obtained in \cite{L}:

\begin{proposition}
In the conditions above there exist a decomposition
$$H_i(T',\lambda'^0) = H'_i(T',\lambda'^0) \bigoplus
H''_i(T',\lambda'^0)$$ into subspaces, invariant under $T''$, such
that:

a/$H'_i(T',\lambda'^0)$ is finite-dimensional (or empty), and the
joint spectrum of operators of $T''$ in $H'_i(T',\lambda'^0)$ (in
the non-empty case) consists on the point $\lambda''^0$.

b/ The joint spectrum of $T''$ in $H''_i(T',\lambda'^0)$ does not
contain the point $\lambda''^0$ .

c/ There exist a finite-dimensional holomorphic subcomplex
$L_*\bc{\lambda'}$ of \newline $K_*\bc{T'-\lambda',X}$ such that
$H_*\bc{L_*\bc{\lambda'^0}}$ coincides with
$H'_i\bc{T'-\lambda'^0,X}$, and
\newline $H_*\bc{K_*\bc{T'-\lambda'^0,X}/L_*\bc{\lambda'^0}}$ - with
$H''_i\bc{T'-\lambda'^0,X}$.

d/ Suppose that $\Pi^{-1}\lambda'^0$ is contained in $\sigma_F(T)$
and therefore consists of finitely many points $\lambda_j =
\bc{\lambda'^0,\lambda''^j}$, $j=1,\ldots,k$. Denote by
$H^j_i(T',\lambda'^0)$ the joint root space of the operators
$T''-\lambda''^j$ in $H_i(T',\lambda'^0)$. Then
$H_i(T',\lambda'^0)=\bigoplus_j H^j_i(T',\lambda'^0)$.
\end{proposition}

In other words, $H'_i(T',\lambda'^0)$ is finite-dimensional and
coincides with the joint root space of the operators
$T_{r+1}-\lambda^0_{r+1},\ldots,T_n-\lambda^0_n$, acting in the
linear space $H_i(T',\lambda'^0)$.

Denote
$$\chi'\bc{T',\lambda'^0}=\sum_{i=0}^r (-1)^i \dim
H'_i(T',\lambda'^0).$$

Let $q$ be the ideal in $\mathcal O_{\lambda^0}$, generated by the
functions $\lambda_1-\lambda^0_1,\ldots , \lambda_r-\lambda^0_r$.
Then

\begin{theorem}

One has $\chi'\bc{T',\lambda'^0}=
e_q\bc{\overline{\mathcal{H}}(T)_{\lambda^0}}$.
\end{theorem}

\begin{proof}

Let, as above, $\mathcal{O}_{\lambda^0}X$ be the stalk at
$\lambda^0$ of the sheaf of germs of $X$-valued holomorphic
functions, and let $\mathcal{O}_{\lambda^0}K_*(T-\lambda,X)$ be
the Koszul complex of the operators of $T-\lambda$ in
$\mathcal{O}_{\lambda^0}X$. This is a complex of
$\mathcal{O}_{\lambda^0}$-modules. In any stage
$\mathcal{O}_{\lambda^0}X_k$ of this complex, $k=0,\ldots,n$, one
can consider the action of the operators $T'-\lambda'^0 =
\bc{T_1-\lambda_1^0,\ldots,T_r-\lambda_r^0}$, and this action
commutes with the differential of the complex. One can form the
Koszul complexes of the operators of $T'-\lambda'^0$ in
$\mathcal{O}_{\lambda^0}X_k$, obtaining a bicomplex of sheaves
with $r+1$ rows and $n+1$ columns. As it is shown below, the
cohomology sheaves of its total complex are in fact supported at
$\lambda^0$ and finite-dimensional, and one may consider its
dimension.

The Euler characteristics of the total complex of the bicomplex
can be computed in two ways. One may consider the homologies
$\mathcal H_i(T)_{\lambda^0}$, $i=0,\ldots,n$ of the initial
complex, and take the alternated sum of dimensions of the
homologies of the Koszul complex of operators $T'-\lambda'^0$ in
it. The alternated sum of these integers for $i=0,\ldots,n$ will
be equal to the alternated sum of the dimensions of the homologies
of the total complex. Since the action of the operators $T_j$ on
$\mathcal H_i(T)_{\lambda^0}$ coincide with the multiplication by
the variable $\lambda_j$, one comes to the Koszul complex of
$\lambda_1-\lambda_1^0,\ldots,\lambda_r-\lambda_r^0$ in $\mathcal
H_i(T)_{\lambda^0}$, and the formula $(*)$ of the p. 1 of the book
\cite{S} can be applied. It shows that for any $i = 0,\ldots,n$
one has
$$\sum_{j=0}^r (-1)^j \dim
H_j\bc{\lambda'-\lambda'^0,{\mathcal H}_i(T)_{\lambda^0}}\ = \
 e_q \bc{\mathcal H_i(T)_{\lambda^0}}$$ and therefore the Euler
characteristic of the total complex is equal to
$e_q\bc{\overline{\mathcal{H}}(T)_{\lambda^0}}$.

On the other hand, take the complex $K_*\bc{T'-\lambda'^0,X}$ with
stages $X_k$ and differentials $d'_k=d'_k\bc{\lambda'^0}: X_k \to
X_{k+1}$, $k=0,\ldots,r$. The Euler characteristic of the total
complex is equal to the alternated sum on $j$, $j=1,\ldots,r$ of
the integers
$$c_j\ = \ \sum_{i=1}^n (-1)^i
\dim {\mathcal H}_i\bc{T-\lambda, {\mathcal
O}_{\lambda^0}H_j\bc{T'-\lambda',X}}$$ where by ${\mathcal
O}_{\lambda^0}H_j\bc{T'-\lambda',X}$ we denote the factor
${\mathcal O}_{\lambda^0}$-module \newline ${\mathcal
O}_{\lambda^0}\ker
d'_j\bc{\lambda'^0}/d'_{j-1}\bc{\lambda'^0}\bc{{\mathcal
O}_{\lambda^0}X_{j-1}}$. Take the finite-dimensional subcomplex
\newline $L_*\bc{\lambda'^0}=\bs{L_j, a_j}$ as in point c/ of prop. 3.2 of \cite{L};
then the action of the operators $T-\lambda$ on the module
${\mathcal O}_{\lambda^0}\bc{X_j/L_j}$ is regular for any $j$, and
therefore in the formula above one can replace ${\mathcal
O}_{\lambda^0}H_j\bc{T'-\lambda',X}$ by ${\mathcal
O}_{\lambda^0}L_j/a_j\bc{{\mathcal O}_{\lambda^0}L_{j-1}}=
{\mathcal O}_{\lambda^0}H'_i\bc{T'-\lambda'^0,X} $. One can easy
prove the following:

\begin{lemma}

Let $H$ be a finite-dimensional space, ${\mathcal O}_{\lambda^0}H$
- the ${\mathcal O}_{\lambda^0}$-module of germs of $H$-valued
holomorphic functions in the point $\lambda^0\in \Cc^n$. Let
$T=\bc{T_1,\ldots ,T_n}$ be an $n$-tuple of operators acting in
$H$ of the form $T_i=\lambda_i^0 + K_i$, $i = 1,\ldots,n$, where
$K_i$ are nilpotent operators in $H$. Then
$$\dim{\mathcal
H}_i\bc{T-\lambda, {\mathcal O}_{\lambda^0}H } \ = \ \left\{
\begin{array}{l}
  0 \quad \mbox{if}\quad i<n \\
  \dim H \quad \mbox{for}\quad i=n. \\
\end{array}
 \right.$$

\end{lemma}

Indeed, if all $K_i$ are zero operators, one obtains the Koszul
complex of the germs $\lambda_i^0 - \lambda_i$ in ${\mathcal
O}_{\lambda^0}H$ exact in all terms except in the last one, where
the module of homology is ${\mathcal
O}_{\lambda^0}H/{\mathfrak{m}}_{\lambda^0}H$, and
${\mathfrak{m}}_{\lambda^0}$ is the maximal ideal in ${\mathcal
O}_{\lambda^0}$. In the general case one can find a filtration of
the linear space $H$ such that all $K_i$ are zero on the
corresponding graduate space, and the assertion follows from the
case of zero operators.

Applying the assertion of the lemma to the space
$H'_i(T',\lambda'^0)$ and the operators
$\lambda_1^0,\ldots,\lambda_r^0, T_{r+1},\ldots,T_n$, one obtains
$c_j = \dim H'_i(T',\lambda'^0)$ for $j=1,\ldots,r$, which proves
the theorem.
\end{proof}

\begin{corollary}
Suppose that $Iso_{max}\bc{\overline{\mathcal{H}}(T)_{\lambda^0}}$
consists of a single prime ideal $p$, (i.e. there is only one
irreducible component $\mathcal P$ of $\sigma_F(T)$ containing
$\lambda^0$), and the coordinates $\lambda' =
\bc{\lambda_1,\ldots,\lambda_r}$ form a coordinate system for
$\sigma_F(T)$ at $\lambda^0$. Then
$$\chi'\bc{T',\lambda'^0}\ =\ e_q\bc{\mathcal O_{\lambda^0}/p}\ l_p(T)\ ,$$
where $q$ is, as above, the ideal generated by
$\lambda_1-\lambda^0_1,\ldots , \lambda_r-\lambda^0_r$.
\end{corollary}

\medskip

\subsection{Computations for the last homology sheaf $\mathcal H_n\bc{T}_{\lambda^0}$}
In this section we will suppose that $\lambda^0 \in \Cc^n$ is such
that the linear subspace of the element of the type $\sum_{i=1}^n
\bc{T_i-\lambda_i^0}x_i$ with $x_1,\ldots,x_n \in X$ is of finite
codimension in $X$. Then, using the arguments of \cite{L}, 4.3.,
it is easy to see that the last homology sheaf $\mathcal
H_n\bc{T}$ of the Koszul complex of $T$ is coherent near
$\lambda^0$.

Let $\frak{m}$ be the maximal ideal in the local ring
$\mathcal{O}_{\lambda^0}$, and denote by $\mathfrak{m}(T)$ the
operator ideal in $L(X)$ generated by the operators
$T_1-\lambda^0_1,\ldots,T_n-\lambda^0_n$. Then the assumption
above means that $\frak{m}(T)X$ is of finite codimension in $X$,
and therefore this will be true for the subspace $\frak{m}^k(T)X$
for any natural $k$.

\begin{theorem}

In the conditions above, one has
$$\dim X/\frak{m}^k(T)X =  \dim \bc { {\mathcal H}_n(T)_{\lambda^0}/\frak{m}^k
{\mathcal H}_n(T)_{\lambda^0}}$$ and therefore for $k$
sufficiently big is a polynomial of $k$ with leading term
\newline $e\bc{\mathcal H_n(T)_{\lambda^0}}\ k^r/r!$, where
$r=\dim \mbox{\rm supp}\, \mathcal H_n(T)_{\lambda^0}$.

\end{theorem}

\begin{proof}

Denote $Y=\frak{m}^k(T)X$ and $Z=X/Y$. Then $Z$ is
finite-dimensional, and the joint spectrum of the operators of $T$
in $Z$ coincides with the point $\lambda^0$. Therefore ${\mathcal
H}_i(T-\lambda,Z)_{\lambda^j}=0$ for all $i<n$. Denote by
${\mathcal H}_{\lambda^0}X$ the ${\mathcal O}_{\lambda^0}$-module
${\mathcal H}_n(T)_{\lambda^0}={\mathcal O}_{\lambda^0}X/\sum_i
\bc{T_i-\lambda_i}{\mathcal O}_{\lambda^j}X$, and by ${\mathcal
H}_{\lambda^0}Y$, ${\mathcal H}_{\lambda^0}Z$ - the corresponding
homology modules for $Y$ and $Z$ respectively. The exact sequence
$0 \to Y \to X \to Z\to 0$ determines an exact sequence of the
corresponding Koszul complexes, and therefore a long exact
sequence of the corresponding homology sheaves, ending with the
sequence $0 \to {\mathcal H}_{\lambda^j}Y\to {\mathcal
H}_{\lambda^j}X \to {\mathcal H}_{\lambda^j}Z \to 0$. We will
prove that the image of ${\mathcal H}_{\lambda^0}Y$ in ${\mathcal
H}_{\lambda^0}X$ coincides with $\frak{m}^k{\mathcal
H}_{\lambda^0}X$. Choose polynomials $g_l(\lambda)$,
$l=1,\ldots,L$, generating the ideal $\frak{m}^k$, and denote by
$\sim$ the relation of equivalence in ${\mathcal O}_{\lambda^0}X$
modulo $\sum \bc{T_i-\lambda_i}{\mathcal O}_{\lambda^0}X$. Any
element of ${\mathcal H}_{\lambda^0}Y$ can be represented by a
germ of $Y$-valued holomorphic function $y(\lambda)=\sum
g_l(T)x_l(\lambda) \in Y$. Then $y$ is equivalent to $\sum
g_l(\lambda)x_l(\lambda)$ and therefore its image in ${\mathcal
H}_{\lambda^0}X$ belongs to $\frak{m}^k {\mathcal
H}_{\lambda^0}X$. We will prove that the image of ${\mathcal
H}_{\lambda^0}Y$ is dense in $\frak{m}^k {\mathcal
H}_{\lambda^0}X$ in the sequential topology (see \cite{G-R},
I.6.2.). Indeed, any element $\xi$ of $\frak{m}^k {\mathcal
H}_{\lambda^0}X$ has the form $\sum_l
g_l(\lambda)\varphi_l(\lambda)$, where $\varphi_l(\lambda)$ are
holomorphic $X$-valued functions, defined in a neighborhood of
$\lambda^0$. $\varphi_l(\lambda)$ can be represented by power
series, converging in some neighborhood of the zero:
$\displaystyle\varphi_l(\lambda) = \sum_{\vec{k}} a_{l, \vec{k}}
\bc{\lambda -\lambda^0}^{\vec{k}}$. Taking the partial sums, one
obtain that $\xi$ is a limit in the sequential topology of the
elements $\displaystyle\xi_N=\sum_l
g_l(\lambda)\sum_{\m{\vec{k}}\leq N} a_{l, \vec{k}} \bc{\lambda
-\lambda^0}^{\vec{k}}$. But the element $\xi_N$ is equivalent to
$\displaystyle\sum_l g_l(T)\sum_{\m{\vec{k}}\leq N}  \bc{T
-\lambda^0}^{\vec{k}}a_{l, \vec{k}}$ and therefore belongs to
${\mathcal H}_{\lambda^0}Y$.

By theorem 10 of \cite{G-R}, II.1.3. all ideals are closed in the
sequential topology, and we obtain that ${\mathcal
H}_{\lambda^0}Y$ and $\frak{m}^k {\mathcal H}_{\lambda^0}X$ are
isomorphic submodules of ${\mathcal H}_{\lambda^0}X$. Taking the
corresponding factor-modules, we obtain an isomorphism between
${\mathcal H}_{\lambda^0}Z$ and ${\mathcal
H}_{\lambda^0}X/\frak{m}^k {\mathcal H}_{\lambda^0}X$. Applying
lemma 2.6, we obtain $$\dim Z = \dim \bc{{\mathcal
H}_n(T)_{\lambda^0}/\frak{m}^k {\mathcal
H}_n(T)_{\lambda^0}}.$$\end{proof}

A similar equality can be proved in a more general situation.
Suppose that, as in th. 2.3,
$f\bc{\lambda}=\bc{f_1\bc{\lambda},\ldots,f_p\bc{\lambda}}$ is a
$k$-tuple of holomorphic functions, defined in the neighborhood of
$\sigma(T)$, $\mu^0 \in \sigma_F(f(T))$, and therefore the
intersection $f^{-1}\bc{\mu^0}\cap \sigma(T)$ is a finite subset
$\bs{\lambda^1,\ldots,\lambda^J}$ of $ \sigma_F(T)$. Let $q$ be
the ideal, generated by
$f_1(\lambda)-\mu^0_1,\ldots,f_p(\lambda)-\mu^0_k$, and $q(T)$ -
the operator ideal in $L(X)$ generated by
$f_1(T)-\mu^0_1,\ldots,f_k(T)-\mu^0_k$. Let $r=\max
\dim_{{\mathcal O}_{\lambda^0}} \bc { {\mathcal
H}_n(T)_{\lambda^j}}$, $j=1,\ldots,J$. Then, applying the theorem
above for the $p$-tuple $f(T)-\mu^0$ and using the functoriality
stated in 2.3, one obtains

\begin{corollary}
In the conditions above, one has
$$\dim X/q^k(T)X = \sum_{j=1}^J \dim\bc { {\mathcal H}_n(T)_{\lambda^j}/q^k {\mathcal
H}_n(T)_{\lambda^j}}$$ and therefore for $k$ sufficiently big is a
polynomial of degree $r$.
\end{corollary}

Recall that the $n$-tuple $T$ has the single value extension
property (SVEP) at $\lambda^0$ if ${\mathcal
H}_i(T)_{\lambda^0}=0$ for all $i \neq n$. Then the theorem
implies

\begin{corollary}

Suppose that, under the conditions of 2.9, the $n$-tuple $T$ has
SVEP at the points $\bs{\lambda^1,\ldots,\lambda^J}$ (or,
equivalently, $f(T)$ has SVEP at $\mu^0$). Then $\dim
\bc{X/q^k(T)X}$ for $m$ sufficiently big is a polynomial of $k$
with leading term $(-1)^n e_q(T) k^r/r!$, where $r$ is the maximum
dimension of the modules ${\mathcal H}_n(T)_{\lambda^j}$ for
$j=1,\ldots,J$.

\end{corollary}

\medskip

\subsection{Application to Arveson`s contraction theory}
Suppose that the operators of $T$ act in the Hilbert space $H$.
Take the polynomials $g_l(\lambda)$, $l=1,\ldots,L$, generating
the ideal $\frak{m}^k$, and denote by $g(T): H^L \to H$ the
operator matrix with entries $g_l(T)$, $l=1,\ldots,L$. Then the
factor-space $H/\frak{m}^k(T)H$ is isomorphic to the kernel of the
operator $g(T) g(T)^* = \sum_{l=1} g_l(T) g_l(T)^*$.

In particular, consider the case when $T$ is a commuting
$n$-contraction with finite rank in the sense of Arveson (see
\cite{A1}), i.e. the operator $1-\sum_{i=1}^n T_i T_i^* $ is
finite-dimensional and positive. Suppose in addition that the
essential spectrum $\sigma_e(T)$ does not contains the origin
(this is obviously satisfied if all the operators $T_i$ are
essentially normal; in this case the essential spectrum is
contained in the unit sphere in $\Cc^n$). Then the Fredholm
spectrum of $T$ is contained in the open unit ball $B_n$. Let
$\phi: \mathcal B(H) \to \mathcal B(H)$ be the Arveson completely
positive map defined by the formula $\phi(A)=\sum_{i=1}^n T_i A
T_i^*$ (see \cite{A2}). In the Arveson`s theory an important role
plays the expression
$$\phi^k(1)\ = \
\sum_{k_1+\ldots+k_n=k} \frac{k!}{ k_1!\ldots k_n!}\
T_1^{k_1}\ldots T_n^{k_n}T_1^{* k_1}\ldots T_n^{* k_n}.$$

Since the functions $\lambda_1^{k_1}\ldots \lambda_n^{k_n}$ with
$k_1+\ldots+k_n=k$ form a system of generators for the $k$-th
power of the maximal ideal in the local ring $\mathcal O_0$ of
germs of holomorphic functions in the point $0\in \Cc^n$, then
from the arguments above it follows:

\begin{proposition}
Let $T$ be a commuting $n$-contraction with finite rank. Then one
has $$\dim \ker \phi^k(1)=\dim H/\frak{m}^k(T)H = \dim \bc {
{\mathcal H}_n(T)_{0}/\frak{m}^k {\mathcal H}_n(T)_{0}}$$ and
therefore $\dim \ker \phi^k(1)$ is a polynomial on $k$ of degree,
equal to the dimension of the ${{\mathcal O}_{0}}$-module
${\mathcal H}_n(T)_{0}$.
\end{proposition}

In \cite{A2}, prop. 7.2, it is proved that $\rank
\bc{1-\phi^k(1)}$ is a polynomial of $k$ of degree $\leq n$. The
degree of this polynomial is called {\em degree of the module} $H$
and denoted by $\deg H$. Since $$\dim \ker \phi^k(1) \leq \rank
\bc{1-\phi^k(1)},$$ then one obtains

\begin{corollary}
Under the assumptions of prop. 2.11 the degree of $H$ is greater
or equal to the dimension of ${\mathcal H}_n(T)_{0}$.
\end{corollary}

\begin{remark}
Suppose that the essential spectrum of $T$ is contained in the
unit sphere in $\Cc^n$. Applying an appropriate M\"{o}bius
transform to the contraction $T$ (see \cite{Gr}), one can replace
in the statement above the point $0$ with an arbitrary point in
the open ball $B_n$. One could make a conjecture that the degree
of a pure finite-rank contraction coincides with the maximal
degree of the corresponding analytic modules, i.e. with the
dimension of the support of ${\mathcal H}_n(T)$; this problem will
be considered elsewhere.
\end{remark}

In  \cite{A2} Arveson introduces the notion of a graded
$n$-contraction as a contraction endowed with a suitable action of
the circle group (for the precise definition see the section 6 of
\cite{A2}). Then, roughly speaking, the Fredholm spectrum and the
corresponding homology sheaves are determined by homogeneous
polynomials in $\Cc^n$. In this case the the inequality of the
proposition above becomes an equality.

\begin{proposition}
Let $T$ be a pure graded finite-rank contraction in the Hilbert
space $H$. Then $$\deg H= \dim_{{\mathcal O}_{0}} {\mathcal
H}_n(T)_{0},$$ and the leading terms of the polynomial $\rank
\bc{1-\phi^k(1)}$ and the Hilbert - Samuel polynomial of
${\mathcal H}_n(T)_{0}$ coincide.
\end{proposition}

\begin{proof}
The assertion  follows almost immediately from the proof of
theorem B of \cite{A2}. Indeed, Arveson constructs a submodule
$H_0$ of finite codimension in $H$ such that $1_{H_0} -
\phi_0^k\bc{1_{H_0}}$ is a (finite-dimensional) projection for any
$k$, $\phi_0^k\bc{A}$ being the Arveson`s completely positive map
for the submodule $H_0$. Therefore $$\rank \bc{1_{H_0} -
\phi_0^k\bc{1_{H_0}}}=\dim \ker \phi_0^k(1_{H_0})$$ and the
equality above is satisfied for the submodule $H_0$.

On the other hand, corollary 1 of theorem C of \cite{A2} shows
that for $k$ sufficiently big the polynomials $\rank \bc{1_{H_0} -
\phi_0^k\bc{1_{H_0}}}$ and $\rank \bc{1_{H} - \phi^k\bc{1_{H}}}$
differ by a polynomial of degree strictly less that $\deg H$.
(Formally, this is stated only in the case when $\deg H=n$, but it
is easy to see that the proof works in the general case as well.)

The same fact is true for the right hand side; indeed, the
arguments used in the proof of 2.8 lead to the exact sequence
$$0 \to {\mathcal H}_n\bc{H_0}_0 \to {\mathcal H}_n\bc{H}_0\to
{\mathcal H}_n\bc{H/H_0}_0\to 0.$$

Since $H/H_0$ is finite-dimensional, then the Hilbert-Samuel
polynomial of \newline ${\mathcal H}_n\bc{H/H_0}_0$ is of degree
zero, and the additivity property (see prop. II.10 of \cite{S})
shows that the Hilbert-Samuel polynomials of ${\mathcal
H}_n\bc{H_0}$ and ${\mathcal H}_n\bc{H}$ differ by a polynomial of
lower degree and therefore have identical leading terms. Now,
since the statement of the proposition is valid for $H_0$, then it
is valid for $H$ also.
\end{proof}

\subsection{Weyl spectrum of $n$-tuple}

\begin{proposition}

Let $\lambda^0 \in \sigma_F(T)$. Then the following assertions are
equivalent:

1/ $z_{\lambda^0}(T)=0$,

2/ $z^{max}_{\lambda}(T)=0$ for $\lambda$ sufficiently close to
$\lambda^0$.

\end{proposition}

\begin{proof}
Suppose that $z_{\lambda^0}(T)=0$, i.e. $l_p\bc{T}=0$ for any
ideal $p \in Iso\bc{\overline{\mathcal H}(T)_{\lambda^0}}$. Then,
as in the proof 2.2, one can see that $l_{\mathcal P}(T) = 0$ for
any irreducible component of $\sigma_F(T)$, containing
$\lambda^0$.

Conversely, suppose that 2/ is satisfied in a neighborhood $U$ of
$\lambda^0$, and $\mathcal P$ is an arbitrary irreducible
component of $\sigma_F(T)$, containing $\lambda^0$. One can choose
a point $\lambda \in \mathcal P \cap U$ such that no other
irreducible component of $\sigma_F(T)$ contains $\lambda$. Then
the condition 2/ in the point $\lambda$ implies that $l_{\mathcal
P}(T) = 0$, and 1/ is satisfied.
\end{proof}

\begin{definition}
The point $\lambda^0 \in \sigma_F(T)$ will be called a {\em Weyl
point} for the $n$-tuple of commuting operators $T$ if the
conditions of the proposition above are fulfilled.

We will denote  by $\rho_\omega(T)$ the set of all Weyl points of
$T$, and by $\omega(T)$ - the {\em Weyl spectrum} of $T$, i.e. the
complement of $\rho_\omega(T)$ in $\sigma(T)$.
\end{definition}

Note that in this definition, unlike in the standard one for a
single operator, the isolated points of the Fredholm spectrum are
not Weyl points.

If the $n$-tuple $T$ possess SVEP, then $z_{\lambda^0}(T)= (-1)^n
z\bc{\mathcal H_n(T)_{\lambda^0}}$ and for any irreducible
components $\mathcal P$ of $\sigma_F(T)$ one has $(-1)^n
l_{\mathcal P}(T)\neq
0$. Therefore, in this case
$\omega(T)=\sigma_e(T)$, and $T$ has "Weyl property".

One can derive a criteria for Weyl points using the spacial
characteristics of the tuple $T$ considered above. Let $\lambda^0
\in \sigma_F(T)$, and suppose that $\mathcal P_1,\ldots,\mathcal
P_s$ are the irreducible components of $\sigma_F(T)$ containing
$\lambda^0$, with dimensions $r_1,\ldots,r_s$ corr. Let $U$ be a
neighborhood of $\lambda^0$ not intersecting other irreducible
components, and take $\lambda^1,\ldots,\lambda^s \in U$ such that
$\lambda^i \in \mathcal P_i$ and $\lambda^i \notin \mathcal P_j$
for $j \neq i$. One can take local coordinate systems in $U$ for
any of components $\mathcal P_i$, $i=1,\ldots,s$, and construct
the corresponding local Euler characteristics
$\chi'\bc{T',\lambda^i}$. Then from 2.75 and 2.15 one obtains:

\begin{proposition}
$\lambda^0$ is a Weyl point iff all $\chi'\bc{T',\lambda^i}=0$,
$i=1,\ldots,s$.
\end{proposition}

The functoriality of $z(T)$, stated in theorem 2.3, implies the
functorial properties of $\omega(T)$.

\begin{proposition}
Let $f\bc{\lambda}=\bc{f_1\bc{\lambda},\ldots,f_k\bc{\lambda}}$ be
a $k$-tuple of holomorphic functions, defined in the neighborhood
of the spectrum $\sigma(T)$ of $T$. Then $\omega\bc{f(T)}\subset
f\bc{\omega(T)}$.
\end{proposition}

\begin{proof}
Take a point $\mu^0 \notin f\bc{\omega(T)}$. Then
$f^{-1}\bc{\mu^0}$, if non-empty, is a finite subset of
$\rho_\omega(T)$, and the same is true for the points $\mu$ in a
sufficiently small neighborhood of $\mu^0$. Then from theorem 2.3
and 2.15 2/ it follows that $\mu^0 \in \rho_\omega\bc{f(T)}$.
\end{proof}

In some particular cases the spectral inclusion in the proposition
above can be replaced by equality; indeed, the following
proposition follows immediately from 2.3 in the same way as above:

\begin{proposition}
Suppose that one of the following two conditions is
satisfied:

a/ $f(z)$ is monomorphic on $\sigma(T)$, or

b/ the operator $T$ has SVEP.

Then $\omega(f(T)) = f(\omega(T))$.
\end{proposition}

\end{document}